\documentclass{article}

\usepackage{amssymb, latexsym}

\usepackage{graphicx}

\usepackage{amsmath}
\usepackage{color}

\newtheorem{theorem}{Theorem}[section]

\newtheorem{corollary}[theorem]{Corollary}

\newtheorem{lemma}[theorem]{Lemma}

\newtheorem{remark}[theorem]{Remark}

\begin{document}
\title{\large\bf  On Harnack inequalities for Witten Laplacian on Riemannian manifolds with super Ricci flows}
\author{\ \ Songzi Li, Xiang-Dong Li
\thanks{Research supported by NSFC No. 11371351, Key Laboratory RCSDS, CAS, No. 2008DP173182, and Hua Luo-Keng Research Grant of AMSS, CAS.} }

\maketitle

\centerline{\it Dedicated to Prof. Ngaiming Mok  for his 60th birthday} 

\vskip1cm

\begin{minipage}{120mm}
{\bf Abstract}.
In this paper, we prove the Li-Yau type Harnack  inequality and Hamilton type dimension free Harnack inequality for the heat equation $\partial_t u=Lu$ associated with the time dependent Witten Laplacian on complete Riemannian manifolds equipped with a variant of the $(K, m)$-super  Perelman Ricci flows and the $K$-super Perelman Ricci flows.

\end{minipage}


\section{Introduction}

\subsection{Motivation}

Differential Harnack inequality is an important topic in the study of heat equations and geometric flows on 
Riemannian manifolds. Let $M$ be an $n$ dimensional complete Riemannian manifold, $u$ be a
positive solution to the heat equation
\begin{eqnarray}
\partial_t u=\Delta u.\label{Heat1}
\end{eqnarray}
In  \cite{LY}, Li and Yau proved that, if the Ricci curvature is bounded from below by a negative constant, i.e., $Ric\geq
-K$, where $K\geq 0$ is a constant, then for all
$\alpha>1$, the following differential  Harnack inequality holds
\begin{eqnarray}
{|\nabla u|^2\over u^2}-\alpha {\partial_t u\over u}\leq
{n\alpha^2\over 2t}+{n\alpha^2K\over \sqrt{2}(\alpha-1)},\ \ \ \forall t>0.\label{LYK}
\end{eqnarray}
In particular, if $Ric\geq 0$, then taking $\alpha\rightarrow 1$,
the Li-Yau Harnack inequality holds
\begin{eqnarray}
{|\nabla u|^2\over u^2}-{\partial_t u\over u}\leq {n\over
2t},\ \ \ \forall t>0.\label{LY}
\end{eqnarray}

In \cite{H1}, on complete Riemannian manifolds with $Ric\geq -K$, Hamilton 
proved  a variant of  the Li-Yau type Harnack inequality  for any positive solution to the heat
equation $(\ref{Heat1})$
\begin{eqnarray}
{|\nabla u|^2\over u^2}-e^{2Kt}{\partial_t u\over u}\leq {n\over
2t}e^{4Kt}, \ \ \ \forall t>0. \label{LYHHar}
\end{eqnarray}
In particular, when $K=0$, the above inequality reduces to the
Li-Yau Harnack inequality $(\ref{LY})$ on complete Riemannian
manifolds with non-negative Ricci curvature.

In the same paper \cite{H1}, Hamilton also proved a dimension free Harnack inequality on
compact Riemannian manifolds with Ricci curvature bounded from
below. More precisely, if $Ric\geq -K$, where $K\geq 0$ is a constant, 
then, for any positive and bounded solution $u$ to the heat equation
$(\ref{Heat1})$, it holds
\begin{eqnarray}
{|\nabla u|^2\over u^2}\leq \left({1\over t}+2K\right)\log\left({A\over u}\right),\ \ \ \forall  t>0, \label{HamHar}
\end{eqnarray}
where $A=\sup\limits\{u(t, x): x\in M, t\geq 0\}$. Indeed, the same result holds for positive solutions (with suitable growth condition)  to the heat equation $\partial_t u=\Delta u$ on complete Riemannian manifolds with
$Ric\geq -K$.

Let $(M, g)$
be a complete Riemannian manifold,
$\phi\in C^2(M)$ (called a potential function), and $d\mu=e^{-\phi}dv$, where $v$ is the Riemannian
volume measure on $(M, g)$. The Witten Laplacian  on $(M, g, \phi)$ is defined  by
\begin{eqnarray*}
L =\Delta -\nabla \phi\cdot\nabla. \label{WL}
\end{eqnarray*}
For all $u, v\in C^\infty_0(M)$, we have the integration by parts formula 
\begin{eqnarray*}
\int_M \langle \nabla u, \nabla v\rangle d\mu=-\int_M Lu vd\mu=-\int_M uLvd\mu.
\end{eqnarray*} By \cite{BE},  for any $u\in C^\infty(M)$, the generalized Bochner formula holds
\begin{eqnarray}
L|\nabla u|^2-2\langle \nabla u, \nabla L u\rangle=2|\nabla^2
u|^2+2Ric(L)(\nabla u, \nabla u), \label{BWF}
\end{eqnarray}
where $\nabla^2 u$ is the Hessian of $u$,  $|\nabla^2 u|$ denotes its Hilbert-Schmidt norm, and 
$$Ric(L)=Ric+\nabla^2\phi.$$
In the literature, $Ric(L)$ is called the (infinite dimensional) Bakry-Emery Ricci curvature associated with the Witten Laplacian $L$ on $(M, g, \phi)$. It plays as a good substitute of the Ricci curvature in many problems in comparison geometry and analysis on complete Riemannian manifolds with smooth weighted volume measures. See \cite{BE, Lot, Li05, Li16} and reference therein.

Following \cite{BE, Li05}, we introduce the $m$-dimensional Bakry-Emery Ricci curvature on $(M, g, \phi)$ by
\begin{eqnarray*}
Ric_{m, n}(L):=Ric+\nabla^2\phi-{\nabla\phi\otimes \nabla \phi\over m-n},
\end{eqnarray*}
where $m\geq n$ is a constant, and $m=n$ if and only if $\phi$ is a constant. When $m=\infty$, we have $Ric_{\infty, n}(L)=Ric(L)$.
Following \cite{BE, Li05}, we say that the Witten Laplacian $L$ satisfies the $CD(K, \infty)$ condition if $Ric(L)\geq K$, and $L$ satisfies the $CD(K, m)$ condition if $Ric_{m, n}(L)\geq K$. Recall that,
when $m\in \mathbb{N}$, the $m$-dimensional Bakry-Emery Ricci curvature $Ric_{m, n}(L)$ has a very natural geometric interpretation.
Indeed, consider the warped product metric on $M^n\times S^{m-n}$ defined by
\begin{eqnarray*}
\widetilde{g}=g_M\bigoplus e^{-{2\phi\over m-n}}g_{S^{m-n}}.\label{WPM}
\end{eqnarray*}
where $S^{m-n}$ is the unit sphere in $\mathbb{R}^{m-n+1}$ with the standard metric $g_{S^{m-n}}$. By \cite{Lot, Li05}, the quantity $Ric_{m, n}(L)$ is equal to the Ricci curvature of the above warped product metric $\widetilde g$ on $M^n\times S^{m-n}$ along the horizontal vector fields.

In \cite{Li05, Li12}, the Li-Yau Harnack inequalities $(\ref{LYK})$ and $(\ref{LY})$ were extended to positive solutions of the heat equation 
\begin{eqnarray}
\partial_t u=Lu, \label{Heat11}
\end{eqnarray}
associated to the 
Witten Laplacian on complete Riemannian manifolds with the $CD(K, m)$-condition for $K\in \mathbb{R}$ and $m\in [n, \infty)$. As application, two-sides Gaussian type heat kernel estimates and the Varadhan short time asymptotic behavior of the heat kernel for the Witten Laplacian were proved in \cite{Li05, Li12}.  
In \cite{Li16},  a slight improved version of Hamilton's Harnack inequality $(\ref{HamHar})$ was proved for any positive and bounded solution to the heat equation $(\ref{Heat11})$ of the 
Witten Laplacian on complete Riemannian manifolds with the $CD(-K, \infty)$-condition. More precisely, letting $(M, g)$ be a complete Riemannian manifold,  $\phi\in C^2(M)$, and assuming that 
$$Ric(L)\geq -K,$$
where $K\geq 0$ is a constant, then for any positive and bounded solution $u$ to the heat equation $(\ref{Heat11})$,  the following optimal dimension free differential Harnack inequality was proved in \cite{Li12}
\begin{eqnarray}
{|\nabla u|^2\over u^2}\leq {2K\over 1-e^{-2Kt}}\log\left({A\over u}\right),\ \ \label{BLH0}
\end{eqnarray}
where $A=\sup\limits\{u(t, x): x\in M, t\geq 0\}$. As far as we know, the above estimate is sharp even for the heat equation $\partial_t u=\Delta u$ on complete Riemannian manifolds with Ricci curvature bounded from below by $-K$, i.e., $Ric\geq -K$. Using the inequality ${2K\over 1-e^{-2Kt}}\leq 2K+{1\over t}$ for $K\geq 0$ and $t>0$, Hamilton's Harnack inequality $(\ref{HamHar})$ for positive and bounded solution to the heat equation $(\ref{Heat11})$ can be derived from $(\ref{BLH0})$.

The aim of this paper is to extend the Li-Yau type and Hamilton type dimension free Harnack
inequalities to positive solutions of the heat equation $(\ref{Heat11})$ for the time dependent Witten Laplacian on complete  Riemannian manifolds equipped with time dependent metrics and potentials. We would like to mention that the Li-Yau and the Hamilton type Harnack inequalities for heat equation $\partial_t u=\Delta u$ on compact or complete Ricci flow has been studied by many authors in the literature. See \cite{Cao, CCLLN, Sun, ZQ2} and references therein. Our aim is to study the Li-Yau  and Harmilton type Harnack inequalities for the heat equation of the time dependent Witten Laplacian on complete Riemannian manifolds equipped with the so-called $K$-super Perelman Ricci flows and a variant of the so-called $(K, m)$-super Ricci flows which we will introduce in Section $1.2$ below.
Indeed, we can also extend the Li-Yau-Hamilton Harnack inequality $(\ref{HamHar})$  to positive solutions of the heat equation $(\ref{Heat11})$ for the time dependent Witten Laplacian on complete  Riemannian manifolds equipped with variant of super Ricci flows. Due to the limit of space, we would like to do this in a separate paper.

\subsection{Statement of main results}

Let $(M, g(t), \phi(t), t\in [0, T])$ be a manifold equipped with a family of time dependent  complete Riemannian metrics $g(t)$ and potential functions $\phi(t)\in C^2(M)$, $t\in [0, T]$. In this paper, we call $(M, g(t), \phi(t), t\in [0, T])$ a $(K, m)$-super Perelman Ricci flow if the metric $g(t)$ and the potential function $\phi(t)$ satisfy the following inequality
\begin{eqnarray*}
{1\over 2}{\partial g\over \partial t}+Ric_{m, n}(L)\geq Kg.
\end{eqnarray*}
When $m=\infty$, i.e., if the metric $g(t)$ and the potential function $\phi(t)$ satisfy the following inequality
\begin{eqnarray*}
{1\over 2}{\partial g\over \partial t}+Ric(L)\geq Kg,
\end{eqnarray*}
we call $(M, g(t), \phi(t), t\in [0, T])$   a $(K, \infty)$-super Perelman Ricci flow or a $K$-super Perelman Ricci flow. 

Note that,  when $\phi(t)\equiv 0$, $t\in [0, T]$, we see that $(M, g(t), \phi(t)\equiv 0, t\in [0, T])$ is a $(K, m)$-super Perelman Ricci flow for any $m\in [n, \infty]$ if and only if 
$(M, g(t), t\in [0, T])$ is the $K$-super Ricci flow in the sense of Hamilton 
\begin{eqnarray*}
{1\over 2}{\partial g\over \partial t}+Ric\geq Kg.
\end{eqnarray*}
On the other hand, we would like to mention that,  the Perelman Ricci flow 
$${1\over 2}{\partial g\over \partial t}+Ric(L)=0$$ 
has been introduced 
in \cite{P1}  as the gradient flow of Perelman's $\mathcal{F}$-functional $\mathcal{F}(g, \phi)=\int_M (R+|\nabla \phi|^2)e^{-\phi}dv$ on $\mathcal{M}\times C^\infty(M)$ under the constraint condition that $d\mu=e^{-\phi}dv$ does not change in time, where $\mathcal{M}$ denotes the set of all Riemannian metrics on $M$. 
 
Our first result, Theorem \ref{Thm0},  extends the dimension free Harnack inequality $(\ref{BLH0})$  to positive and bounded solutions of the heat equation 
$\partial_t u=Lu$  for time dependent Witten Laplacian on manifolds equipped with a  complete $(-K)$-super Perelman Ricci flow. As far as we know,  our result is new even in 
the case of super Ricci flow without potential, i.e., $\phi(t)=0$, $t\in [0, T]$. See Theorem \ref{Thm1}.

\begin{theorem}\label{Thm0} Let $M$ be a manifold equipped with a family of time dependent  complete Riemannian metrics and $C^2$-potentials
$(g(t), \phi(t), t\in [0, T])$. Suppose that $(M, g(t), \phi(t), t\in [0, T])$ is a $(-K)$-super Perelman
Ricci flow 
\begin{eqnarray}
{1\over 2}{\partial g\over \partial t}+Ric(L)\geq -Kg, \label{SRFK}
\end{eqnarray}
where $K\geq 0$ is a constant independent of  $t\in [0, T]$. Let $u$ be a positive and bounded solution to the heat equation
\begin{eqnarray*}
\partial_t u=Lu,
\end{eqnarray*}
where
$$L=\Delta_{g(t)}-\nabla_{g(t)}\phi(t)\cdot\nabla_{g(t)}
$$
is the time dependent Witten Laplacian on $(M, g(t), \phi(t))$.
Suppose that 
\begin{eqnarray}
\int_0^T \int_M  \left( \left|\nabla\left({|\nabla u|^2\over u}\right)(y)\right|^2+|\nabla (u\log u)|^2(y) \right)p_{0, t}(x, y)d\mu(y)dt<\infty,\label{energy}
\end{eqnarray}
where $p_{s, t}(x, y)$ denotes the fundamental solution to the heat equation $\partial_t u=Lu$ with respect to the weighted volume measure $\mu$,  $0\leq s\leq t\leq T$. 
Then for all $x\in M$ and $t>0$,
\begin{eqnarray}
{|\nabla u|^2\over u^2}\leq {2K\over 1-e^{-2Kt}}\log\left({A\over u}\right),\ \ \label{BLH}
\end{eqnarray}
where $A=\sup\limits\{u(t, x): x\in M, t\geq 0\}$.  Using the inequality ${2K\over 1-e^{-2Kt}}\leq 2K+{1\over t}$ for $K\geq 0$ and $t>0$, we have the Hamilton Harnack inequality 
\begin{eqnarray}
{|\nabla u|^2\over u^2}\leq \left({1\over
t}+2K\right)\log\left({A\over u}\right).\label{Ham}
\end{eqnarray}
In particular, when $K=0$, i.e.,
$(M, g(t), \phi(t), t\in [0, T])$ is a  manifold equipped with a complete super Perelman
Ricci flow
\begin{eqnarray*}
{1\over 2}{\partial g\over \partial t}+Ric(L)\geq 0,
\end{eqnarray*}
we have
\begin{eqnarray*}
{|\nabla u|^2\over u^2}\leq {1\over t}\log\left({A\over u}\right).\ \ \label{HH1}
\end{eqnarray*}
\end{theorem}

In particular, when $\phi(t)\equiv 0$, $t\in [0, T]$, we have the following dimension free Harnack inequality for the heat equation $(\ref{Heat1})$ of the time dependent Laplace-Beltrami on manifolds with $(-K)$-super Ricci flows.  

\begin{theorem}\label{Thm1} Let $(M, g(t), t\in [0, T])$ be a  manifold equipped with a complete $(-K)$-super 
Ricci flow 
\begin{eqnarray*}
{1\over 2}{\partial g\over \partial t}+Ric\geq -Kg, \label{SRFK2}
\end{eqnarray*}
where $K\geq 0$ is a constant independent of  $t\in [0, T]$. Let $u$ be a positive and bounded solution to the heat equation associated with the time dependent Laplace-Beltrami 
\begin{eqnarray*}
\partial_t u=\Delta u.
\end{eqnarray*}
Suppose that 
\begin{eqnarray*}
\int_0^T \int_M  \left( \left|\nabla\left({|\nabla u|^2\over u}\right)(y)\right|^2+|\nabla (u\log u)|^2(y)|  \right)p_{0, t}(x, y)dvol(y)dt<\infty,
\end{eqnarray*}
where $p_{s, t}(x, y)$ denotes the fundamental solution to the heat equation $\partial_t u=\Delta u$ with respect to the volume measure $dvol$, $0\leq s\leq t\leq T$.  
Then for all $x\in M$ and $t>0$,
\begin{eqnarray*}
{|\nabla u|^2\over u^2}\leq {2K\over 1-e^{-2Kt}}\log\left({A\over u}\right),\ \ \label{BLH2}
\end{eqnarray*}
where $A=\sup\limits\{u(t, x): x\in M, t\geq 0\}$.  Using the inequality ${2K\over 1-e^{-2Kt}}\leq 2K+{1\over t}$ for $K\geq 0$ and $t>0$, we have the Hamilton Harnack inequality 
\begin{eqnarray*}
{|\nabla u|^2\over u^2}\leq \left({1\over
t}+2K\right)\log\left({A\over u}\right).\label{Ham2}
\end{eqnarray*}
In particular, when $K=0$, i.e.,
$(M, g(t), t\in [0, T])$ is a  manifold equipped with a complete  super Ricci flow
\begin{eqnarray*}
{1\over 2}{\partial g\over \partial t}+Ric\geq 0,
\end{eqnarray*}
we have
\begin{eqnarray}
{|\nabla u|^2\over u^2}\leq {1\over t}\log\left({A\over u}\right).\ \ \label{HH2}
\end{eqnarray}
\end{theorem}

Integrating the differential Harnack inequality $(\ref{BLH})$ along the geodesic on  $(M, g(t))$, we have the following Harnack inequality for positive solutions of the heat equation of the time dependent Witten Laplacian on super Perelman Ricci flows.

\begin{corollary}\label{cor1} Under the same condition and notation as in Theorem \ref{Thm0} and Theorem \ref{Thm1}, for any $\delta>0$, and for all $x, y\in M$, $0<t<T$, we have
\begin{eqnarray*}
u(x, t)\leq A^{\delta\over 1+\delta} u(y, t)^{1\over 1+\delta}\exp\left\{{1+\delta^{-1}\over 4(1+\delta)}{2K\over 1-e^{-2Kt}}d_t^2(x, y)\right\}.
\end{eqnarray*}
where $d_t(x, y)$ denotes the distance between $x$ and $y$ in $(M, g(t))$. 
In particular, when $K=0$, we have
\begin{eqnarray*}
u(x, t)\leq A^{\delta\over 1+\delta} u(y, t)^{1\over 1+\delta}\exp\left\{{1+\delta^{-1}\over 4(1+\delta)}{d_t^2(x, y)\over t}\right\}.
\end{eqnarray*}
\end{corollary}

The next result extends the Li-Yau type Harnack inequality to positive solutions of the heat equation $\partial_t u=Lu$ for time dependent Witten Laplacian on compact or 
complete Riemannian manifolds equipped with  a variant of the $(K, m)$-super Perelman Ricci flows. 

\begin{theorem}\label{LYHSRF-A}  Let $(M, g(t), t\in [0, T])$ be a Riemannian manifold with a family of time dependent complete Riemannian metrics $g(t)$ and potentials $\phi(t)\in C^2(M)$, $t\in [0, T]$. Let $L=\Delta_{g(t)}-\nabla_{g(t)}\phi(t)\cdot\nabla_{g(t)}$, and $u$ be a positive solution to the heat equation $\partial_t u=Lu$. Let $\partial_t g=2h$ and $\alpha>1$. Suppose that there exist two constants $K\geq 0$ and $m>n$ independent of $t\in [0, T]$  such that
\begin{eqnarray}
{1\over 2}(1-\alpha)\partial_t g+Ric_{m, n}(L)\geq -Kg,\label{mmm1}
\end{eqnarray}
and assume that $A^2=\max\limits \left[ |h|^2+{({\rm Tr}h)^2\over m-n}\right]<\infty$ and $B=\max\limits |S|<\infty$, where
\begin{eqnarray*}
S(\cdot)=2h(\nabla \phi, \cdot\rangle-\langle 2{\rm div} h-\nabla {\rm Tr}_g h+\nabla \partial_t \phi, \cdot\rangle+{2 {\rm Tr} h\over m-n}\langle \nabla\phi, \cdot\rangle.
\end{eqnarray*}
If $M$ is  compact, then for any $\gamma>0$ and  $t\in (0, T]$, we have

\begin{eqnarray*}
{|\nabla u|^2\over u^2}-\alpha {\partial_t u\over u}\leq {m\alpha^2\over 4t}\left[1+\sqrt{1+{t^2\over m}\left(4A^2+{m(2K+\gamma)^2\over (\alpha-1)^2}+{2B^2\over \gamma}\right)}\right].
\end{eqnarray*}
If $(M, g(t), t\in [0, T])$ are complete, then for any $\gamma>0$ and $t\in (0, T]$, we have
\begin{eqnarray}
{|\nabla u|^2\over u^2}-\alpha {\partial_t u\over u}\leq {m\alpha^2\over
4t}\left[1+C_4(K_2+\sqrt{K_1})t+\sqrt{(1+C_4(K_2+\sqrt{K_1})t)^2+{Dt^2\over m}}\right].\ \ \ \  \label{Harcomplete-L}
\end{eqnarray}
where $C_4$ is a constant depending only on $m$,  $K_1$ and $K_2$ are two positive constants such that $Ric_{m, n}(L)\geq -K_1$ and $h\geq -K_2$, and $D=4A^2+{m(2K+\gamma)^2\over (\alpha-1)^2}+{B^2
\over 2\gamma}$.
\end{theorem}

When $\phi(t)\equiv 0$, $t\in [0, T]$, we have the following Li-Yau Harnack inequality  for positive solutions of the heat equation $\partial_t u=\Delta_{g(t)} u$ on complete Riemannian manifolds equipped with  
a variant of the $K$-super Ricci flows.
 
\begin{theorem}\label{LYHSRF-B}  Let $(M, g(t), t\in [0, T])$ be a manifold equipped with a family of time dependent complete Riemannian metrics $g(t)$.  Let $u$ be a positive solution to the heat equation 
$$\partial_t u=\Delta_{g(t)} u.$$ 
Let $\partial_t g=2h$ and $\alpha>1$. Suppose that there exist two constants $K\geq 0$ and $m>n$ independent of $t\in [0, T]$  such that
\begin{eqnarray}
{1\over 2}(1-\alpha)\partial_t g+Ric\geq -Kg,\label{mmmm1}
\end{eqnarray}
and assume that $A^2=\max\limits \left[ |h|^2+{({\rm Tr}h)^2\over m-n}\right]<\infty$ and $B=\max\limits |S|<\infty$, where
\begin{eqnarray*}
S(\cdot)=-\langle 2{\rm div} h-\nabla {\rm Tr}_g h, \cdot\rangle.
\end{eqnarray*}
If $M$ is  compact, then for any $\gamma>0$ and $t\in (0, T]$, we have
\begin{eqnarray}
{|\nabla u|^2\over u^2}-\alpha {\partial_t u\over u}\leq {m\alpha^2\over 4t}\left[1+\sqrt{1+{t^2\over m}\left(4A^2+{m(2K+\gamma)^2\over (\alpha-1)^2}+{2B^2\over \gamma}\right)}\right]. \label{HarnCom}
\end{eqnarray}
If $(M, g(t), t\in [0, T])$ are complete, then for any $\gamma>0$ and $t\in (0, T]$, we have
\begin{eqnarray}
{|\nabla u|^2\over u^2}-\alpha {\partial_t u\over u}\leq {m\alpha^2\over
4t}\left[1+C_4(m)(K_2+\sqrt{K_1})t+\sqrt{(1+C_4(m)(K_2+\sqrt{K_1})t)^2+{Dt^2\over m}}\right].  \label{Harcomplete}
\end{eqnarray}
where $C_4$ is a constant depending only on $m$,  $K_1$ and $K_2$ are two positive constants such that $Ric\geq -K_1$ and $h\geq -K_2$, and $D=4A^2+{m(2K+\gamma)^2\over (\alpha-1)^2}+{B^2
\over 2\gamma}$.
\end{theorem}

By standard method as in Li-Yau \cite{LY} and Chow et al \cite{CCLLN}, integrating the above Li-Yau differential Harnack quantity along paths on the space-time, we can derive the following parabolic Harnack inequality for the solution of the heat equation on different points in space-time.

\begin{corollary}\label{cor2}  Let $(M, \widetilde{g})$ be a complete Riemannian manifold, $(g(t), \phi(t), t\in [0, T])$ be a family of complete Riemannian metrics and $C^2$-potentials on $M$. Assuming that for each $t\in [0, T]$, 
there exists a constant $C>0$ such that
\begin{eqnarray*}
C^{-1}\widetilde{g}\leq g(t)\leq C\widetilde{g}.
\end{eqnarray*} Let $u$ be a positive solution to the heat equation $\partial_t u=Lu$. Then, under the same condition and notation as in Theorem \ref{LYHSRF-A} or Theorem \ref{LYHSRF-B}, for any $\alpha>1$, $x_1, x_2\in M$ and $0<t_1<t_2\leq T$, we have
\begin{eqnarray*}
{u(x_2, t_2)\over u(x_1, t_1)}\leq e^{-C_7 (t_2-t_1)}\left({t_1\over t_2}\right)^{m\alpha\over 2}\exp\left(-{C\alpha\over 4 }{d_{\widetilde{g}}^2(x_1, x_)\over t_2-t_1}\right), 
\end{eqnarray*}
where $C_7=C_4(K_2+\sqrt{K_1})+\sqrt{D\over m}$ and $D={(2K+\gamma)^2\over (\alpha-1)^2}+{2B^2
\over m\gamma}+{4A^2\over m}$.
\end{corollary}

To end this Section, let us give some remarks and compare our results with known results in the literature.

\begin{remark}{\rm \ \ 

\begin{itemize}

\item In \cite{ZQ2}, Q. Zhang proved $(\ref{HH2})$ for positive and bounded solutions to the heat equation $\partial_t u=\Delta_{g(t)} u$ on complete Riemmanian
manifolds equipped with the Ricci flow $\partial_t g=-2Ric$. By a probabilistic approach,  Guo, Philliposki and Thalmaier  \cite{Guo} proved the Hamilton Harnack inequality $(\ref{HH2})$ for
  the backward heat equation $\partial_t u=-\Delta_{g(t)} u$ on Riemannian manifolds equipped with complete backward super Ricci flow $\partial_t g\leq 2Ric$.  See also \cite{Cao} for  Hamilton type Harnack inequality on Riemannian manifolds equipped with  complete Ricci flow with $|Ric|\leq K$. In \cite{Li16}, the second author give two proofs of the optimal Hamilton dimension free Harnack inequality $(\ref{BLH0})$ for positive and bounded solution to the heat equation $\partial_t u=Lu$
 on complete Riemannian manifolds $Ric(L)\geq -K$. The proof of Theorem \ref{Thm0} is similar to the probabilistic proof of $(\ref{BLH0})$ given in \cite{Li16}. 
  For another proof of  Theorem \ref{Thm0} derived from the reversal logarithmic Sobolev inequality on complete $K$-super Perelman Ricci flow, see \cite{LL16}. 
 
\item A local version of the Li-Yau Harnack inequality in Theorem \ref{LYHSRF-A} and Theorem \ref{LYHSRF-B} is proved in Section 4.2, see Theorem \ref{LYHlocal}. 

\item In \cite{Sun}, J. Sun proved the Li-Yau Harnack inequality for positive solutions of the heat equation $\partial_t u=\Delta_{g(t)}u$ on manifold $M$  with a family of complete Riemannian 
metrics $(g(t), t\in [0, T])$. The assumptions in \cite{Sun}  are given by: $\partial_t g=2h$, $Ric\geq -K_1g$, $-K_2g\leq h\leq K_3g$, and $|\nabla h|\leq K_4$.  Under these conditions, for any $\alpha>1$, Sun proved that 
\begin{eqnarray*}
{|\nabla u|^2\over u^2}-\alpha {\partial_t u\over u}\leq {n\alpha^2\over
t}+C(K_1+K_2+K_3+K_4+\sqrt{2K_4}),
\end{eqnarray*}
where $C$ depends only on $n$ and $\alpha$. Indeed, under his conditions, we have 
\begin{eqnarray*}
{1\over 2}(1-\alpha)\partial_t g+Ric\geq -(K_1+(\alpha-1)K_3)g,
\end{eqnarray*}
i.e., $(\ref{mmmm1})$ holds with $K=K_1+(\alpha-1)K_3$. Moreover, 
\begin{eqnarray*}
|h|^2\leq n(K_2+K_3)^2.
\end{eqnarray*}
Using the inequality $|{\rm Tr} h|^2\leq n |h|^2$, we have
\begin{eqnarray*}
A^2=\max\limits \left[ |h|^2+{({\rm Tr}h)^2\over m-n}\right]  \leq  {mn\over m-n} (K_2+K_3)^2.
\end{eqnarray*}
On the other hand
\begin{eqnarray*}
B=|2{\rm div} h-
\nabla {\rm Tr}_g h|=|2g^{ij}\nabla_{i}h_{jl}-g^{ij}\nabla_l h_{ij}|\leq 3|g||\nabla h|\leq 3\sqrt{n}K_4.
\end{eqnarray*} 
Therefore, Theorem \ref{LYHSRF-B} applies and yields the following estimate
\begin{eqnarray*}
{|\nabla u|^2\over u^2}-\alpha {\partial_t u\over u}&\leq& {m\alpha^2\over 
4t}\left[1+C_4(K_2+\sqrt{K_1})t\right]\\
& &+ {m\alpha^2\over 4t}\sqrt{(1+C_4(K_2+\sqrt{K_1})t)^2+{t^2\over m}\left(4A^2+{m(2K+\gamma)^2\over (\alpha-1)^2}+{B^2
\over 2\gamma}\right)} \\
&\leq& {m\alpha^2\over 2t}+C_{m, n, \alpha, \gamma}(1+\sqrt{K_1}+K_1+K_2+K_3+K_4),
\end{eqnarray*}
where $C_{m, n, \alpha, \gamma}$ depends only on $m, n, \alpha$ and $\gamma$. Taking $m=2n$, we have
\begin{eqnarray*}
{|\nabla u|^2\over u^2}-\alpha {\partial_t u\over u}\leq {n\alpha^2\over t}+C_{n, \alpha, \gamma}(1+\sqrt{K_1}+K_1+K_2+K_3+K_4).
\end{eqnarray*}

\item In the case of the Ricci flow, i.e., $\partial_t g=-2Ric$, we have $h=-Ric$, ${\rm Tr}h=-R$, where $R$ is the scalar curvature.  In this case, ${1\over 2}(1-\alpha)\partial_t g+Ric=\alpha Ric$. Thus $(\ref{mmmm1})$ reads as 
\begin{eqnarray*}
\alpha Ric\geq -Kg.
\end{eqnarray*}
Note that the second Bianchi identity says that
\begin{eqnarray*}
{\rm div} Ric-{1\over 2}\nabla R=0.
\end{eqnarray*}
Thus $B\equiv 0$ for all $t\in [0, T]$.  Let $u$ be a positive solution to the heat equation $\partial_t u=\Delta_{g(t)}u$ on a Ricci flow $(M, g(t), t\in [0, T])$  with $Ric\geq -\alpha^{-1}Kg$. In 
the case $M$ is  compact, then for any $t\in (0, T]$, we have
\begin{eqnarray}
{|\nabla u|^2\over u^2}-\alpha {\partial_t u\over u}\leq {m\alpha^2\over 4t}\left[1+\sqrt{1+{4t^2\over m}\left(A^2+{mK^2\over (\alpha-1)^2}\right)}\right]. 
\label{compactHar}
\end{eqnarray}
and in the case  $(M, g(t), t\in [0, T])$ are complete, then for any  $t\in (0, T]$, we have
\begin{eqnarray}
{|\nabla u|^2\over u^2}-\alpha {\partial_t u\over u}&\leq &{m\alpha^2\over
4t}\left[1+C_4(K_2+\sqrt{K_1})t\right]\nonumber\\
& &+{m\alpha^2\over
4t}\sqrt{(1+C_4(K_2+\sqrt{K_1})t)^2+{4t^2\over m}\left(A^2+{mK^2\over (\alpha-1)^2}\right)}, \ \ \ \ \label{completeHar}
\end{eqnarray}
where $C_4$ is a constant depending only on $m$, $K_1\geq 0$ and $K_2\geq 0$ are two constant such that $Ric\geq -K_1g$ and $h\geq -K_2g$, i.e., $-K_1g\leq Ric\leq{ K_2\over 2}g$. 
See also Bailesteanu-Cao-Pulemotov \cite{Cao} for the Li-Yau type Harnack estimate on complete Riemannian manifolds with Ricci flow such that $|Ric|\leq K$.

\item In the case of the backward Ricci flow, i.e., $\alpha=2$, $\partial_t g=2Ric$, we have $h=Ric$, ${\rm Tr}h=R$, where $R$ is the scalar curvature.  In this case,  ${1\over 2}(1-\alpha)\partial_t g+Ric=(2-\alpha) Ric$, and $(\ref{mmmm1})$ reads as 
\begin{eqnarray*}
(2-\alpha) Ric\geq -Kg.
\end{eqnarray*}
By the second Bianchi identity,  $B\equiv 0$ for all $t\in [0, T]$.  Let $u$ be a positive solution to the heat equation $\partial_t u=\Delta_{g(t)}u$ on a backward Ricci flow $(M, g(t), t\in [0, T])$  
with $(2-\alpha)Ric\geq -Kg$.  Then, if $M$ is  compact, for any $t\in (0, T]$, $(\ref{compactHar})$ holds, and if $(M, g(t), t\in [0, T])$ are complete, $(\ref{completeHar})$ holds. 
In particular, if $(M, g(t), t\in [0, T])$ is a backward Ricci flow and $\alpha=2$, we have $K=0$. In this case,   if $M$ is compact, then for any  $t\in (0, T]$,  we have
\begin{eqnarray*}
{|\nabla u|^2\over u^2}-\alpha {\partial_t u\over u}\leq {m\alpha^2 \over 4t}\left[1+\sqrt{1+{4A^2t^2\over m}}\right],
\end{eqnarray*}
and if $(M, g(t), t\in [0, T])$ are complete, then  for any  $t\in (0, T]$,  we have
\begin{eqnarray}
{|\nabla u|^2\over u^2}-\alpha {\partial_t u\over u}\leq {m\alpha^2\over
4t}\left[1+C_4(K_2+\sqrt{K_1})t+\sqrt{(1+C_4(K_2+\sqrt{K_1})t)^2+{4A^2t^2\over m}}\right],\label{kkk}
\end{eqnarray}
where $C_4$ is a constant depending only on $m$, and $K_1\geq 0$ and $K_2\geq 0$ are two constant such that $Ric\geq -K_1g$ and $h\geq -K_2g$,  i.e., $Ric\geq -{K_2\over 2}g$.
Thus $K_2=2K_1$ and $(\ref{kkk})$ reads as follows
\begin{eqnarray*}
{|\nabla u|^2\over u^2}-\alpha {\partial_t u\over u}\leq {m\alpha^2\over
4t}\left[1+C_4(K_1+\sqrt{K_1})t+\sqrt{(1+C_4(K_1+\sqrt{K_1})t)^2+{4A^2t^2\over m}}\right].\label{kkkk}
\end{eqnarray*}

\item In the case $g(t)$ and $\phi(t)$ are independent of $t\in [0, T]$, we have $A=B=0$, and $K_2=0$.  Thus, on any compact or complete Riemannian manifold with  $Ric_{m, n}(L)\geq -Kg$, for all $t\in (0, T]$, we have
\begin{eqnarray*}
{|\nabla u|^2\over u^2}-\alpha {\partial_t u\over u}\leq {m\alpha^2\over 4t}\left[1+C_4\sqrt{K}t+\sqrt{(1+C_4\sqrt{K}t)^2+{4K^2t^2\over (\alpha-1)^2}}\right].
\end{eqnarray*}
Hence  
\begin{eqnarray}
{|\nabla u|^2\over u^2}-\alpha {\partial_t u\over u}\leq {m\alpha^2\over 2t}\left[1+C_4\sqrt{K}t+{Kt\over \alpha-1}\right]. \label{LYHmK}
\end{eqnarray}
From the proof of Theorem \ref{LYHSRF-B}, we see that $C_4=C(m-1)$ for some constant $C>0$. 
In particular, on any  complete Riemannian manifold with $Ric_{m, n}(L)\geq 0$, we recapture the generalized Li-Yau Harnack inequality for any positive solution to the heat equation $\partial_t u=Lu$ (see \cite{Li05, Li12})
\begin{eqnarray*}
{|\nabla u|^2\over u^2}-{\partial_t u\over u}\leq {m\over 2t}. \label{LYHm0}
\end{eqnarray*}
Taking $m=n$, $\phi\equiv 0$, and $L=\Delta$, then, for any positive solution to the heat equation $\partial_t u=\Delta u$ on a complete Riemannian manifold with $Ric\geq -K$, and 
for any $\alpha>1$, we have (compare to $(\ref{LYK})$)
\begin{eqnarray*}
{|\nabla u|^2\over u^2}-\alpha {\partial_t u\over u}\leq {n\alpha^2\over 2t}\left[1+C_4\sqrt{K}t+{Kt\over \alpha-1}\right], \label{LYHnK}
\end{eqnarray*}
Here $C_4=C(n-1)$ for some constant $C>0$. 
In particular, we recapture the Li-Yau Harnack inequality $(\ref{LY})$ for any positive solution to the heat equation $\partial_t u=\Delta u$ on any complete Riemannian manifold with non-negative Ricci curvature.
\end{itemize}
}
\end{remark}

The rest of this paper is organized as follows.  In Section $2$, we prove Theorem \ref{Thm0} and Corollary \ref{cor1}.  In Section $3$,  
we prove the Li-Yau Harnack inequality for the heat equation of time dependent Witten Laplacian on compact Riemannian manifolds equipped with 
a variant of the $(K, m)$-super Perelman Ricci flows, i.e., Theorem \ref{LYHSRF-B}.  
In Section $4$,  we prove Theorem \ref{LYHSRF-B}  on complete Riemannian manifolds equipped with a variant of the $(K, m)$-super Perelman Ricci flows.  

This paper is an improved version of a part of the authors' preprint \cite{LL16} which will be divided into three papers due to the limit of the space.  See also \cite{LL17a, LL17b}.

\section{Hamilton type  Harnack inequality on super Ricci flows}

To prove the Hamilton type Harnack inequality on  $(-K)$-super Perelman Ricci flows, we extend the probabilistic approach which was used in time independent case in \cite{Li16}. First we introduce the $L$-diffusion process  on $(M, g(t), t\in [0, T])$. Following \cite{ACT}, let $(U_t, t\in [0, T])$ be the solution of the Stratonovich SDE on the orthonormal frame bundle $(O(M), g(t), t\in [0, T])$ over $(M, g(t), t\in [0, T])$
\begin{eqnarray*}
dU_t&=&\sum\limits_{i=1}^n H_i (U_t)\circ dW_t^i-\left[(\nabla \phi)^{H}(U_t)+\sum\limits_{\alpha, \beta=1}^n{\partial g\over \partial t}(U_te_\alpha, U_te_\beta)V_{\alpha, \beta}(U_t)\right]dt,\\
U_0&=&u\in (O(M), g_0),
\end{eqnarray*}
where $\{H_i\}_{i=1}^n$ denote the canonical vector fields on $(O(M), g(t))$, $\{V_{\alpha,\beta}\}_{\alpha, \beta=1}^n$ denote the canonical vertical vector fields on $(O(M), g(t))$, and $(\nabla\phi)^H$ denotes the horizontal lift of the vector field $\nabla \phi$ from $(M, g(t))$ to $(O(M), g(t))$. The $L$-diffusion process on $(M, g(t))$ is defined by \begin{eqnarray*}
X_t=\pi (U_t)
\end{eqnarray*}
By \cite{ACT}, for smooth $f\in [0, T]\times M\rightarrow \mathbb{R}$, Ito's formula holds
\begin{eqnarray*}
df(t, X_t)=\left(\partial_t +L\right)f(t, X_t)dt+\sum\limits_{i=1}^n \nabla_{e_i} f(t, X_t)dW_t^i.
\end{eqnarray*}

\noindent{\bf Proof of Theorem \ref{Thm0}}. By direct calculation and the generalized Bochner formula, we have
\begin{eqnarray}
(\partial_t-L){|\nabla u|^2\over u}=-{2\over u}\left|\nabla^2
u-{\nabla u\otimes \nabla u\over u}\right|^2-{2\over u }\left({1\over 2}{\partial
g\over
\partial t}+Ric(L)\right)(\nabla u, \nabla u).\label{Boh1}
\end{eqnarray}
Thus, on manifold with a $(-K)$-super Perelman Ricci flow,  we have 
\begin{eqnarray*}
(\partial_t-L){|\nabla u|^2\over u}\leq 2K{|\nabla u|^2\over u}.
\end{eqnarray*}
Note that
\begin{eqnarray*}
(\partial_t-L) (u\log (A/u))={|\nabla u|^2\over u}.
\end{eqnarray*}
Let $\psi(t)={1-e^{-2Kt}\over 2K}$. Then $
\psi'(t)+2K\psi(t)=1$. Define
\begin{eqnarray*}
h(x, t):=\psi(t){|\nabla u|^2\over u}-u\log (A/u).
\end{eqnarray*}
Then at $t=0$, $h\leq 0$,  and for $t>0$, it holds
\begin{eqnarray*}
(\partial_t-L)h\leq \left[\psi'(t)+2K\psi(t)-1\right]{|\nabla u|^2\over u}= 0.
\end{eqnarray*}

In the case $M$ is compact, the maximum principle yields that $h(x,
t)\leq 0$ for all time $t>0$ and $x\in M$. In the case $(M, g(t), t\in [0, T])$ is a
complete non-compact Riemannian manifold with a $(-K)$-super Perelman Ricci flow, we
can give a probabilistic proof to $(\ref{BLH})$ as follows. Let $X_t$ be the $L$-diffusion
process on $(M, g(t))$ starting from $X_0=x$. Applying It\^o's formula to
$h(X_t, T-t)$, $t\in [0, T]$, we have
\begin{eqnarray*}
h(X_t, T-t)=h(X_0, T)+\int_0^t \nabla h(X_s, T-s)dW_s+\int_0^t \left(L-{\partial\over \partial t}\right)h(X_s, T-s)ds,
\end{eqnarray*}
where the second term in the right hand side is the It\^o's stochastic integral with respect to the Brownian motion $\{W_s, s\in [0, t]\}$.
In particular, taking $t=T$, we obtain
\begin{eqnarray*}
h(X_T, 0)=h(X_0, T)+\int_0^T \nabla h(X_s, T-s)\cdot dW_s+\int_0^T \left(L-{\partial\over \partial t}\right)h(X_s, T-s)ds.
\end{eqnarray*}
Note that, under the condition $(\ref{energy})$, we have 
\begin{eqnarray*}
\mathbb{E}\left[\int_0^T |\nabla h(X_s, T-s)|^2 ds\right]<\infty.
\end{eqnarray*}
Hence $M_t= \int_0^t \nabla h(X_s, T-s)dW_s$ is a martingale with respect to $\mathcal{F}_t=\sigma(W_s, s\leq t)$. 
Taking the expectation on both sides, the martingale property of It\^o's integral implies that
\begin{eqnarray*}
E\left[h(X_T, 0)\right]=h(x, T)+E\left[\int_0^T \left(L-{\partial\over \partial t}\right)h(X_s, T-s)ds\right]\geq h(x, T).
\end{eqnarray*}
As $h(y, 0)\leq 0$ for all $y\in M$, we derive that $h(x, T)\leq 0$ for all $T>0$ and $x\in M$.  \hfill $\square$

\medskip

\noindent{\it Proof of Corollary \ref{cor1}}. The proof is similar to the one of Theorem $1.1$ in \cite{H1}.   Let $l(x, t)=\log{A/u(x, t)}$. Then the differential Harnack  inequality $(\ref{BLH})$ in Theorem \ref{Thm0} implies
\begin{eqnarray*}
|\nabla\sqrt{l(x, t)}|={1\over 2}{|\nabla l(x, t)|\over \sqrt{l(x, t)}}\leq {1\over 2}\sqrt{2K\over 1-e^{-2Kt}}.
\end{eqnarray*}
Fix $x, y\in M$ and integrate along a geodesic on $(M, g(t))$ linking $x$ and $y$, the above inequality yields
\begin{eqnarray*}
\sqrt{\log{A/u(x, t)}}\leq \sqrt{\log{A/u(y, t)}}+{1\over 2}\sqrt{2K\over 1-e^{-2Kt}}d_t(x, y).
\end{eqnarray*}
where $d_t(x, y)$ denotes the distance between $x$ and $y$ in $(M, g(t))$.
Combining this with the elementary inequality
\begin{eqnarray*}
(a+b)^2\leq (1+\delta)a^2+(1+\delta^{-1})b^2,
\end{eqnarray*}
we can derive the desired Harnack inequality for $u$ in  Corollary \ref{cor1}. \hfill $\square$

\section{Li-Yau Harnack inequality on compact super Perelman Ricci flows}

Let $u$ be a positive solution to the heat equation $\partial_t u=Lu$. Let $f=\log u$. Then
\begin{eqnarray*}
(L-\partial_t) f=-|\nabla f|^2.
\end{eqnarray*}
Let
\begin{eqnarray*}
F=t(|\nabla f|^2-\alpha f_t).
\end{eqnarray*}

\subsection{The commutator $[\partial_t,  L]f$}

Let $M$ be a  manifold with a family of time dependent metrics $(g(t), t\in [0, T])$ and potentials $\phi(t)\in C^2(M)$, $t\in [0, T]$. Let $\partial_t g=2h$.

\begin{lemma} \label{LLLL} For any $f\in C^\infty(M)$, we have
\begin{eqnarray*}
\partial_t |\nabla f|^2=-2h(\nabla f, \nabla f)+2\langle \nabla f, \nabla f_t\rangle,
\end{eqnarray*}
and
\begin{eqnarray}
[\partial_t, L] f=-2\langle h, \nabla^2 f\rangle+2h(\nabla \phi, \nabla f)-\langle 2 {\rm div}h-\nabla {\rm Tr}_g h+\nabla \partial_t \phi, \nabla f\rangle. \label{commu}
\end{eqnarray}
\end{lemma}
{\it Proof}. By direct calculation, we have
\begin{eqnarray*}
\partial_t |\nabla f|^2=\partial_t g^{ij}(t)\nabla_i f\nabla_j f=\partial_t g^{ij}(t)\nabla_i f\nabla_j f+2g^{ij}(t)\nabla_i f\nabla_j f_t.
\end{eqnarray*}
Note that
\begin{eqnarray*}
\partial_t g^{ij}(t)=-\partial_t g_{ij}(t)=-2h_{ij}.
\end{eqnarray*}
The first equality follows. On the other hand, by \cite{CLN, Sun}, we have
\begin{eqnarray*}
\partial_t \Delta_{g(t)} f=\Delta_{g(t)}\partial_t f-2\langle h, \nabla^2 f\rangle-2\langle {\rm div}h- {1\over 2}\nabla{\rm Tr}_g h, \nabla f\rangle.
\end{eqnarray*}
Combining this with
\begin{eqnarray*}
\partial_t \langle \nabla \phi, \nabla f\rangle=-\partial_t g(\nabla \phi, \nabla f)+\langle \nabla\phi_t, \nabla f\rangle+\langle \nabla \phi, \nabla f_t\rangle,
\end{eqnarray*}
we obtain  $(\ref{commu})$ in Lemma \ref{LLLL}.  \hfill $\square$
%

\begin{lemma}\label{lemm}
\begin{eqnarray}
(L-\partial_t) F&=&2t\left( |\nabla^2 f|^2+(Ric(L)+(1-\alpha)h)(\nabla f, \nabla f)\right)\nonumber\\
& &\ \ \ \ \ \ \ \ -2\langle \nabla f, \nabla F\rangle-t^{-1}F+\alpha t[\partial_t, L]f.    \label{lemm}
\end{eqnarray}
\end{lemma}
{\it Proof}. By the Bochner formula and using
\begin{eqnarray*}
\partial_t |\nabla f|_{g(t)}^2=-\partial_t g(t)(\nabla f, \nabla f)+2\langle \nabla f, \nabla f_t\rangle_{g(t)},
\end{eqnarray*}
we have
\begin{eqnarray*}
LF&=&tL|\nabla f|^2-\alpha t Lf_t\\
&=&2t\left( |\nabla^2 f|^2+Ric(L)(\nabla f, \nabla f)+\langle \nabla f, \nabla L f\rangle \right) -\alpha t L\partial_t f\\
&=&2t\left( |\nabla^2 f|^2+Ric(L)(\nabla f, \nabla f)+\langle \nabla f, \nabla (f_t-|\nabla f|^2)\rangle \right) -\alpha t L\partial_t f\\
&=&2t\left( |\nabla^2 f|^2+Ric(L)(\nabla f, \nabla f)\right)-2\langle \nabla f, \nabla F\rangle+2(1-\alpha) t\langle \nabla f, \nabla f_t\rangle-\alpha t L\partial_t f\\
&=&2t\left( |\nabla^2 f|^2+Ric(L)(\nabla f, \nabla f)\right)-2\langle \nabla f, \nabla F\rangle\\
& &\hskip1cm +2t(1-\alpha) h(\nabla f, \nabla f)+(1-\alpha)t \partial_t |\nabla f|^2-\alpha t L\partial_t f.
\end{eqnarray*}
On the other hand
\begin{eqnarray*}
\partial_t F&=&(|\nabla f|^2-\alpha f_t)+t\partial_t |\nabla f|^2-\alpha t f_{tt}\\
&=&(|\nabla f|^2-\alpha f_t)+t\partial_t |\nabla f|^2-\alpha t \partial_t (Lf+|\nabla f|^2)\\
&=&(|\nabla f|^2-\alpha f_t)+(1-\alpha)t\partial_t |\nabla f|^2-\alpha t \partial_t Lf.
\end{eqnarray*}
Combing above formulas, we derive $(\ref{lemm})$. \hfill $\square$

\medskip

\begin{lemma}\label{lemmaF} For any $\alpha>1$, we have
\begin{eqnarray}
(L-\partial_t) F&\geq &{2F^2\over \alpha^2 m t}+{4(\alpha-1)\over
m\alpha^2}|\nabla f|^2F
+{2t(\alpha-1)^2\over m\alpha^2}|\nabla f|^4-{t\alpha^2\over 2}\left[{({\rm Tr}h)^2\over m-n}+|h|^2 \right]\nonumber\\
& &+2t(Ric_{m, n}(L)+(1-\alpha)h)(\nabla f, \nabla f)-2\langle
\nabla f, \nabla F\rangle-t^{-1}F+\alpha t S(\nabla
f).\nonumber \\
& &\hfill \label{ooooo}
\end{eqnarray}
\end{lemma}
{\it Proof}. Substituting $[\partial_t, L]f$ into $(\ref{lemm})$, we have
\begin{eqnarray*}
(L-\partial_t) F&=&2t\left|\nabla^2 f-{\alpha h\over 2}\right|^2-{t\alpha^2|h|^2\over 2} +2t(Ric(L)+(1-\alpha)h)(\nabla f, \nabla f)\\
& &-2\langle \nabla f, \nabla F\rangle-t^{-1}F+\alpha t S_1(\nabla f),
\end{eqnarray*}
where
\begin{eqnarray*}
S_1(\nabla f)=2h(\nabla \phi, \nabla f\rangle-\langle 2{\rm div} h-\nabla {\rm Tr}_g h+\nabla \phi_t, \nabla f\rangle.
\end{eqnarray*}
Using the inequality $|S|^2\geq {1\over n}|{\rm Tr} S|^2$ for $n\times n$  symmetric matrices $S$ and the  Cauchy-Schwartz inequality $(a+b)^2\geq {a^2\over 1+\varepsilon}-{b^2\over \varepsilon}$ for all $\varepsilon>0$, we can obtain
\begin{eqnarray*}
\left|\nabla^2 f-{\alpha h\over 2}\right|^2&\geq &{1\over n}\left|\Delta f-{\alpha {\rm Tr}h\over 2} \right|^2\\
&\geq& {|Lf|^2\over n(1+\varepsilon)}-{\left|\nabla \phi \cdot\nabla f-{\alpha {\rm Tr}h \over 2}\right|^2\over n\varepsilon}.
\end{eqnarray*}
Let $m:=n(1+\varepsilon)$. Then
\begin{eqnarray}
(L-\partial_t) F&\geq &{2t\over m}|Lf|^2-{2t\over m-n}\left|\nabla \phi \cdot\nabla f-{\alpha {\rm Tr}h \over 2}\right|^2 -{t\alpha^2|h|^2\over 2} +2t(Ric(L)+(1-\alpha)h)(\nabla f, \nabla f)\nonumber\\
& &-2\langle \nabla f, \nabla F\rangle-t^{-1}F+\alpha t S_1(\nabla f)\nonumber\\
&=&{2t\over m}|Lf|^2-{t\alpha^2 ({\rm Tr}h)^2\over 2(m-n)}-{t\alpha^2|h|^2\over 2} +2t(Ric_{m, n}(L)+(1-\alpha)h)(\nabla f, \nabla f)\nonumber\\
& &-2\langle \nabla f, \nabla F\rangle-t^{-1}F+\alpha t S_1(\nabla f)+{2\alpha t {\rm Tr} h\over m-n}\langle \nabla\phi, \nabla f\rangle.\label{FF1}
\end{eqnarray}
Let
\begin{eqnarray*}
S(\cdot)=S_1(\cdot)+{2 {\rm Tr} h\over m-n}\langle \nabla\phi, \cdot\rangle.
\end{eqnarray*}
Substituting $Lf=|\nabla f|^2-f_t={F\over \alpha t}+{\alpha-1\over \alpha}|\nabla f|^2$ into $(\ref{FF1})$, we have
\begin{eqnarray*}
(L-\partial_t) F
&\geq &{2t\over m}\left[{F\over \alpha t}+{\alpha-1\over \alpha}|\nabla f|^2\right]^2-{t\alpha^2\over 2}\left[{({\rm Tr}h)^2\over m-n}+|h|^2 \right]\\
& &+2t(Ric_{m, n}(L)+(1-\alpha)h)(\nabla f, \nabla f)-2\langle \nabla f, \nabla F\rangle-t^{-1}F+\alpha t S(\nabla f).
\end{eqnarray*}
This completes the proof of Lemma \ref{lemmaF}.  \hfill $\square$

\medskip

Note that $A^2=\max\limits \left[ |h|^2+{({\rm Tr}h)^2\over m-n}\right]$, $B=\max\limits |S|$. Under the assumption $(\ref{mmm1})$,  i.e.,  $Ric_{m, n}(L)+(1-\alpha)h\geq -K$, we have
\begin{eqnarray*}
(L-\partial_t) F&\geq &{2F^2\over \alpha^2 m t}+{2t(\alpha-1)^2\over m\alpha^2}|\nabla f|^4-{t\alpha^2A^2\over 2}\\
& &-2Kt|\nabla f|^2-2\langle \nabla f, \nabla F\rangle-t^{-1}F-\alpha B t |\nabla f|.
\end{eqnarray*}
Using the inequality
\begin{eqnarray*}
ax^4+bx^2+cx\geq -{(b-\gamma)^2\over 4a}-{c^2\over 4\gamma}, \label{xx1}
\end{eqnarray*}
where $\gamma>0$ is any positive constant, we can derive that
\begin{eqnarray*}
{2t(\alpha-1)^2\over m\alpha^2}|\nabla f|^4-2Kt|\nabla f|^2-\alpha B
t |\nabla f|\geq -{m\alpha^2 t(2K+\gamma)^2\over
8(\alpha-1)^2}-{\alpha^2 B^2 t\over 4\gamma}.\label{0}
\end{eqnarray*}
Hence
\begin{eqnarray}
(L - \partial_t)F
\geq  {2F^2\over \alpha^2 mt}-\frac{F}{t} - 2\langle \nabla f, \nabla F\rangle -{t\alpha^2A^2\over 2}-{m\alpha^2 t(2K+\gamma)^2\over 8(\alpha-1)^2} -{\alpha^2 B^2 t\over 4\gamma}. \label{llllF}
\end{eqnarray}

\subsection{Proof of Theorem \ref{LYHSRF-A} in compact case}

When $M$ is compact. Let $(x_0, t_0)$ be the point where $F$ achieves the maximum on $M\times [0, T]$.
Then $\nabla F(x_0, t_0)=0$, $\Delta F(x_0, t_0)\leq 0$ and $\partial_t F(x_0, t_0)\geq 0$. Therefore, at $(x_0, t_0)$,
\begin{eqnarray*}
(L-\partial_t)F\leq 0.
\end{eqnarray*}
By $(\ref{llllF})$, we have
\begin{eqnarray*}
0\geq {2F^2\over \alpha^2 m t_0}-{F\over t_0} -{t_0\alpha^2A^2\over 2}-{m\alpha^2 t_0(2K+\gamma)^2\over 8(\alpha-1)^2}-{\alpha^2 B^2 t_0\over 4\gamma}.
\label{xx2}
\end{eqnarray*}
This yields, for any  $t\in (0, T]$, we have
\begin{eqnarray*}
F&\leq& {m\alpha^2\over 4}\left[1+\sqrt{1+{t_0^2\over m}\left(4A^2+{m(2K+\gamma)^2\over (\alpha-1)^2}+{2B^2\over \gamma}\right)}\right]\\
&\leq&{m\alpha^2\over 4}\left[1+\sqrt{1+{T^2\over m}\left(4A^2+{m(2K+\gamma)^2\over (\alpha-1)^2}+{2B^2\over \gamma}\right)}\right].
\end{eqnarray*}
In particular, at time $t=T$, we derive the Li-Yau Harnack inequality in Theorem \ref{LYHSRF-B}
\begin{eqnarray*}
{|\nabla u|^2\over u^2}-\alpha {\partial_t u\over u} \leq {m\alpha^2\over
4t}\left[1+\sqrt{1+{t^2\over m}\left(4A^2+{m(2K+\gamma)^2\over (\alpha-1)^2}+{2B^2\over \gamma}\right)}\right].
\end{eqnarray*}
 \hfill $\square$

\section{Proof of Theorem \ref{LYHSRF-A} in complete case}

\subsection{A lemma} 

Fix $o\in M$. Let $Q_{2R, T}=\{(x, t)\in M\times [0, T]: d(x, o, t)\leq 2R, t\in [0, T]\}$. Let $\eta\in C^2([0, \infty), [0, 1]) $ be such that $\eta(r)=1$ on $[0, 1]$, $\eta=0$ on $[2, \infty)$, $0\leq \eta \leq 1$ on $[1, 2]$, $\eta'(r)\leq 0$, $|\eta'(r)|^2\leq C_1\eta(r)$ and $\eta''(r)\geq -C_2$, where $C_1$ and $C_2$ are two positive constants. Define
\begin{eqnarray*}
\psi(x, t)=\psi(d(x, o, t))=\eta\left({d(x, o, t)\over R}\right)=\eta\left({\rho(x, t)\over R}\right),
\end{eqnarray*}
where $\rho(x, t)=d(x, o, t)$ denotes the geodesic distance between $x$ and $o$ on $(M, g(t))$.

We need the Laplacian comparison theorem on manifolds with time dependent metrics and potentials.

\begin{lemma}\label{comparison2} Let $M$ be a complete Riemannian manifold equipped with a family of time dependent complete Riemannian metrics $g(t)$ and potentials $\phi(t)$, $t\in [0, T]$. Let
$$\partial_t g=2h.$$
Suppose that
 $Ric_{m, n}(L)\geq -K_1$, $h\geq -K_2$, where $K_1, K_2$ are  two positive constants. Then
\begin{eqnarray*}
(L-\partial_t)\psi\geq -C_1K_2\psi^{1/2}-{C_1\over R}(m-1)\sqrt{K_1}\coth(\sqrt{K_1}\rho)-{C_2\over R^2}.
\end{eqnarray*}
\end{lemma}
{\it Proof}. By \cite{Li05},  as $Ric_{m, n}(L)\geq -K_1$, the following Laplacian comparison theorem holds
$$
Ld(x_0, x, t)\leq (m-1)\sqrt{K_1}\rho\coth(\sqrt{K_1}\rho),
$$
and
\begin{eqnarray*}
L\psi&=&\eta'(d(x_0, x, t)/R){L d(x_0, x, t)\over R}+\eta''(d(x_0, x, t)/R){|\nabla d(x_0, x, t)|^2\over R^2}\\
&\geq &-{C_1\over R}(m-1)\sqrt{K_1}\coth(\sqrt{K_1}\rho)-{C_2\over R^2}.
\end{eqnarray*}
On the other hand, let $\gamma: [a, b]\rightarrow M$ be a fixed path such that $\gamma(a)=x$ and $\gamma(b)=y$. Let $S=\dot \gamma(s)$. Given a time $t_0\in [0, T]$, assuming that $\gamma$ is parameterize by the arc length with respect to metric $g(t_0)$ on $M$, then $|S|=1$ at time $t=t_0$. Moreover,  the evolution of the length of $\gamma$ with respect to $g(t)$ is given by
\begin{eqnarray*}
\left.{d\over dt}\right|_{t=t_0}L_{g(t)}(\gamma)&=&\int_a^b \left.{d\over dt}\right|_{t=t_0}\sqrt{g(t)(S, S)}ds\\
&=&{1\over 2}\int_a^b \left.{\partial_tg(t)(S, S)\over \sqrt{g(t)(S, S)}}\right|_{t=t_0}ds\\
&=&{1\over 2}\int_a^b \left.{\partial g(t)\over \partial t}(S, S)\right|_{t=t_0}ds.
\end{eqnarray*}
This yields, under the assumption $h\geq -K_2$, where $K_2\geq 0$,
\begin{eqnarray*}
\partial_t d(x, y, t)= \int_a^b h(S, S)ds\geq -K_2d(x, y, t).
\end{eqnarray*}
Since $-C_1\eta^{1/2}(r)\leq \eta'(r)\leq 0$, and $K_2\geq 0$, it holds
\begin{eqnarray*}
-\partial_t \psi&=&-{\eta'(\rho/R)\partial_t d(x_0, x, t)\over R}\\
&\geq&{\eta'(\rho/R) K_2 d(x_0, x, t)\over R}\\
&\geq &-{C_1K_2\over R}\psi^{1/2} d(x_0, x,t).
\end{eqnarray*}
Combining this with the lower bound of $L\psi$, we have
\begin{eqnarray*}
(L-\partial_t)\psi\geq -C_1K_2\psi^{1/2}-{C_1\over R}(m-1)\sqrt{K_1}\coth(\sqrt{K_1}\rho)-{C_2\over R^2}.
\end{eqnarray*}
The proof of Lemma \ref{comparison2} is completed. \hfill $\square$

\subsection{The local version of the Li-Yau Harnack inequality}

In this subsection we prove a local version of the Li-Yau type Harnack inequality for positive solutions to the heat equation $\partial_t u=Lu$ on complete Riemannian manicolds with a variant of the $(K, m)$-super Perelman Ricci flows.  More precisely, we have the following
\begin{theorem}\label{LYHlocal} Let $(M, g(t), \phi(t), t\in [0, T])$ be a manifold equipped with a family of time dependent complete Riemannian metrics and $C^2$-potentials. Under the same condition as  in Theorem \ref{LYHSRF-B},  for any $\alpha>1$ and $R>0$, we have the following local Li-Yau type differential  Harnack inequality on  $Q_{2R, T}=\{(x, t)\in M\times [0, T]: d(x, o, t)\leq 2R, t\in [0, T]\}$
\begin{eqnarray}
{|\nabla u|^2\over u^2}-\alpha {\partial_t u\over u} &\leq& {m\alpha^2\over 4t}\left[1+Et+\sqrt{\left(1+Et \right)^2+{Dt^2\over m} }\right],  \label{local}
\end{eqnarray}
where $E=C_4(K_2+\sqrt{K_1})+{C_5\over R}+{C_6\over
R^2}$,  $D=4A^2+{m(2K+\gamma)^2\over (\alpha-1)^2}+{2B^2\over \gamma}$,  $C_4$, $C_5$ are constants depending only on $m$, and $C_6$ is a constant 
depending only on $m$ and $\alpha$.
\end{theorem}
{\it Proof}.  Let $F=t(|\nabla \log u|^2-\alpha \partial_t \log u)$. Since $\rho$ is Lipschitz on the complement of the cut locus of $o$,
$\psi$ is a Lipschitz function with support in $Q_{2R, T}$. As
explained in Li and Yau [43], an argument of Calabi \cite{Cal58}
allows us to apply the maximum principle to $\psi F$. Let $(x_0,
t_0)\in M\times [0, T]$ be a point where $\psi F$ achieves the
maximum. Then, at $(x_0, t_0)$,
\begin{eqnarray*}
\partial_t(\psi F)\geq 0, \ \Delta(\psi F)\leq 0, \ \nabla(\psi F)=0,
\end{eqnarray*}
which yields
\begin{eqnarray*}
(L-\partial_t)(\psi F)=\Delta (\psi F)-\nabla\phi\cdot\nabla(\psi
F)-\partial_t(\psi F)\leq 0.
\end{eqnarray*}Note that
\begin{eqnarray*}
(L-\partial_t)(\psi F)=\psi(L-\partial_t)F+(L-\partial_t)\psi
F+2\nabla\psi \cdot\nabla F.
\end{eqnarray*}
By Lemma \ref{comparison2}, we have
\begin{eqnarray*}
(L-\partial_t)\psi\geq -C_1K_2\psi^{1/2}-{C_1\over
R}(m-1)\sqrt{K_1}\coth(\sqrt{K_1}\rho)-{C_2\over R^2}.
\end{eqnarray*}
Therefore, at $(x_0, t_0)$, we have
\begin{eqnarray}
0\geq \psi(L-\partial_t)F+2\nabla\psi \cdot\nabla F-A(R,
T)F,\label{mmm}
\end{eqnarray}
where
$$A(R, T):=C_1K_2\psi^{1/2}+{C_1\over R}(m-1)\sqrt{K_1}\coth(\sqrt{K_1}\rho)+{C_2\over R^2}.$$
Denote $$C_3=A(R,
T)+2|\nabla\psi|^2\psi^{-1}.$$ Using $\sqrt{K_1}\rho\coth(\sqrt{K_1}\rho)\leq 1+\sqrt{K_1}\rho$, we have
\begin{eqnarray*}
C_3&\leq& C_1K_2+{C_1(m-1)(1+\sqrt{K_1}R)\over R}+{2C_1+C_2\over R^2}\\
&\leq&C_1(K_2+(m-1)\sqrt{K_1})+{C_1(m-1)\over R}+{2C_1+C_2\over R^2}.
\end{eqnarray*}
To simplify the notation, write
\begin{eqnarray*}
C_3\leq C_4(K_2+\sqrt{K_1})+{C_5\over R}+{C_6\over R^2}.
\end{eqnarray*}
Note that, at $(x_0, t_0)$, $\nabla\psi\cdot\nabla F=-\psi |\nabla
\psi|^2 F$. Substituting $(\ref{ooooo})$ into $(\ref{mmm})$, at
$(x_0, t_0)$, we have
\begin{eqnarray*}
0&\geq&\psi(L-\partial_t)F-A(R, T)F+2\nabla\psi\cdot\nabla F\\
&\geq &\psi(L-\partial_t)F -(A(R, T) + 2|\nabla\psi|^2\psi^{-1})F\\
&\geq& {2\psi F^2\over m\alpha^2 t}-\left({\psi \over
t}+C_3\right)F+{4(\alpha-1) \psi |\nabla f|^2 F\over
m\alpha^2}-{2C_2\over R}\psi^{1/2}|\nabla f|F\\
& &+\psi t\left[{2(\alpha-1)^2\over m\alpha^2}|\nabla f|^4-2K|\nabla
f|^2-\alpha B |\nabla f|-{\alpha^2 A^2\over 2}\right].
\end{eqnarray*}
By the inequality $ax^2-bx\geq {4b^2\over a}$ and $(\ref{llllF})$, and
multiplying the both sides by $\psi t_0$, we have
\begin{eqnarray*}
0&\geq& {2(\psi F)^2\over m\alpha^2}-\left(\psi +C_3t+{m\alpha^2 C_2^2 t\over 4(\alpha-1)R^2}\right)\psi F\\
& &-{\alpha^2\psi^2 t^2\over 8}\left(4A^2+{m(2K+\gamma)^2\over (\alpha-1)^2}+{2B^2\over \gamma}\right).
\end{eqnarray*}
Let $D=4A^2+{m(2K+\gamma)^2\over (\alpha-1)^2}+{2B^2\over \gamma}$. We see that, at any $(x, t) \in Q_{R, T}$, we have
\begin{eqnarray*}
F(x,t) &\leq& (\psi F)(x_0, t_0)\\
&\leq&{m\alpha^2\over 4}\left[1+C_3t_0+{m\alpha^2 C_2^2 t_0\over
4(\alpha-1)R^2}+\sqrt{\left(1+C_3t_0+{m\alpha^2 C_2^2 t_0\over
4(\alpha-1)R^2} \right)^2+{D\psi^2 t_0^2\over m}}\right]\\
&\leq&{m\alpha^2\over 4}\left[1+C_3T+{m\alpha^2 C_2^2 T\over
4(\alpha-1)R^2}+\sqrt{\left(1+C_3T+{m\alpha^2 C_2^2 T\over
4(\alpha-1)R^2} \right)^2+{D\psi^2 T^2\over m}}\right].
\end{eqnarray*}
In particular, taking $t=T$, we have 
\begin{eqnarray*}
F(x,t) &\leq& {m\alpha^2\over 4}\left[1+\left(C_4(K_2+\sqrt{K_1})+{C_5\over R}+{C_6\over
R^2}+{m\alpha^2 C_2^2 \over
4(\alpha-1)R^2}\right)t\right]\\
& &+ {m\alpha^2\over 4}\sqrt{\left(1+\left(C_4(K_2+\sqrt{K_1})+{C_5\over R}+{C_6\over
R^2}+{m\alpha^2 C_2^2 \over 4(\alpha-1)R^2}\right)t \right)^2+{Dt^2\over m} }.
\end{eqnarray*}
This completes the proof of the  Li-Yau Harnack inequality on $Q_{2R, T}$. \hfill $\square$

\subsection{Proof of Theorem \ref{LYHSRF-A} and Corollary \ref{cor2}}

\noindent{\it Proof of Theorem \ref{LYHSRF-A}}.\ When $(M, g(t), t\in [0, T])$ are  complete non-compact, taking $R\rightarrow \infty$ in $(\ref{local})$, we obtain the Li-Yau differential  Harnack inequality on $M\times (0, T]$
\begin{eqnarray*}
{|\nabla u|^2\over u^2}-\alpha {\partial_t u\over u} \leq {m\alpha^2\over
4t}\left[1+C_4(K_2+\sqrt{K_1})t+\sqrt{(1+C_4(K_2+\sqrt{K_1})t)^2+{Dt^2\over m}} \right].
\end{eqnarray*}
This completes the proof of  Theorem \ref{LYHSRF-A}. 
 \hfill $\square$

\bigskip

\noindent{\it Proof of Corollary \ref{cor2}}. Let $\gamma: [t_1, t_2]\rightarrow M$ be a smooth path with $\gamma(t_i)=x_i$, $i=1, 2$.  Then
\begin{eqnarray*}
\log {u(x_2, t_2)\over u(x_1, t_1)}&=&\int_{t_1}^{t_2} {d\over dt}\log u(\gamma(t), t)dt\\
&=&\int_{t_1}^{t_2} \left( \partial_t \log u +\langle \nabla \log u, \dot \gamma(t)\rangle\right) dt.
\end{eqnarray*}
By the Li-Yau Harnack differential inequality in Theorem \ref{LYHSRF-A}, we have
\begin{eqnarray*}
\log {u(x_2, t_2)\over u(x_1, t_1)}&\geq& \int_{t_1}^{t_2} \left({1\over \alpha} |\nabla \log u|^2-{m\alpha \over 4t}\left(2(1+(C_4(K_2+\sqrt{K_1})+\sqrt{Dm^{-1}})t\right)+\langle \nabla \log u, \dot \gamma(t)\rangle\right)dt\\
&\geq& -{\alpha \over 4}\int_{t_1}^{t_2} |\dot \gamma(t)|^2dt-{m\alpha\over 2}\log\left({t_2\over t_1}\right)-C_7(t_2-t_1).
\end{eqnarray*}
Therefore, for any path $\gamma$ on $M$ with $\gamma(t_i)=x_i$, $i=1, 2$, we have
\begin{eqnarray*}
{u(x_2, t_2)\over u(x_1, t_1)}\geq  e^{-C_7(t_2-t_1)} \left( {t_1\over t_2} \right)^{m\alpha\over 2} \exp \left(-{\alpha\over 4} \int_{t_1}^{t_2} |\dot \gamma(t)|_{g(t)}^2dt\right).
\end{eqnarray*}
Let $\gamma(t)$ be a constant speed minimal geodesic linking $x_1$ and $x_2$ on $(M, \widetilde{g})$. By assumption, we have 
\begin{eqnarray*}
	\left| \dot \gamma(t) \right|^2_{g(t)} \leq C\left| \dot \gamma(t) \right|^2_{\widetilde{g}} = {C d^2_{\widetilde{g}} (x_1, x_2)\over (t_2-t_1)^2}.
\end{eqnarray*}
This yields
\begin{eqnarray*}
{u(x_2, t_2)\over u(x_1, t_1)}\geq  e^{-C_7(t_2-t_1)} \left( {t_1\over t_2} \right)^{m\alpha\over 2} \exp \left(-{C\alpha\over 4} {d^2_{\widetilde{g}}(x_1, x_2)\over t_2-t_1}\right).
\end{eqnarray*}
\hfill $\square$

\bigskip

\noindent{\bf Acknowledgement}.  The authors would like to thank Prof. A. Thalmaier and Dr. Yuzhao Wang for useful discussions and their interests on this paper.

\begin{flushleft}

Songzi Li, School of Mathematical Science, Beijing Normal University, No.~19, Xin Jie Kou Wai Da Jie, 100875, China, Email: songzi.li@bnu.edu.cn
\medskip

Xiang-Dong Li, Academy of Mathematics and Systems Science, Chinese
Academy of Sciences, 55, Zhongguancun East Road, Beijing, 100190, China, 
E-mail: xdli@amt.ac.cn

and 

School of Mathematical Sciences, University of Chinese Academy of Sciences, Beijing, 100049, China

\end{flushleft}

\end{document}